\documentclass[aima]{amsart}

\usepackage{amssymb}
\usepackage{amsmath}
\usepackage{cite}

\usepackage{mathrsfs}
\usepackage{pstricks,pst-node}
\psset{linewidth=0.4pt}  

\title{The representation type of Hecke algebras of type $B$}
\author{Susumu Ariki}
\address{Research Institute for Mathematical Sciences\\
Kyoto University, Kyoto 606-8502. Japan.\\
E--mail: ariki@kurims.kyoto-u.ac.jp 
}
\author{Andrew Mathas}
\address{School of Mathematics and Statistics F07\\
University of Sydney, Sydney NSW 2006. Australia.\\
E--mail: a.mathas@maths.usyd.edu.au\\
www: www.maths.usyd.edu.au/u/mathas
}

\def\And{\text{\ and\ }}
\def\For{\text{\ for\ }}
\def\ForAll{\text{\ for all\ }}
\def\If{\text{\ if\ }}
\def\Where{\text{\ where\ }}

{\catcode`\|=\active
  \gdef\set#1{\mathinner{\lbrace\,{\mathcode`\|"8000%
                                   \let|\midvert #1}\,\rbrace}}
}
\def\midvert{\egroup\mid\bgroup}

\makeatletter
\def\Number#1{\refstepcounter{equation}
              \leqno(\theequation)\if*#1%
              \else\def\@currentlabel{{\rm\theequation}}\label{#1}%
              \fi}
\makeatother

\newtheorem{Theorem}[equation]{Theorem}
\newtheorem{Proposition}[equation]{Proposition}
\newtheorem{Lemma}[equation]{Lemma}
\newtheorem{Corollary}[equation]{Corollary}

\newenvironment{Point}[1]%
  {\ifx*#1\let\pointlabel\relax\else\def\pointlabel{#1}\fi
   \refstepcounter{equation}\trivlist
   \item[\hskip\labelsep\bf\theequation
         \ifx\pointlabel\relax\else\space\pointlabel\space\fi]
   \ignorespaces\it
  }{\relax}

\def\Sum{\displaystyle\sum}

\def\){\big)}
\def\({\big(}
\let\>\rangle
\let\<\langle
\let\incl\hookrightarrow
\let\0\varnothing              

\let\bar\overline

\let\gedom\trianglerighteq
\let\gdom\vartriangleright

\def\End{\mathop{\rm End}\nolimits}
\def\Hom{\mathop{\rm Hom}\nolimits}

\def\N{{\mathbb N}}
\def\Z{{\mathbb Z}}
\def\Q{{\mathbb Q}}

\let\To\longrightarrow
\def\map#1#2{\,{:}\,#1\!\longrightarrow\!#2}

\newcounter{CaseCounter}
\numberwithin{CaseCounter}{equation}
\def\Case#1{\addtocounter{CaseCounter}{1}\smallskip
            \noindent{\it Case \arabic{CaseCounter}:\space#1.\space}}

\def\A{\mathscr A}

\def\H{\mathscr H}
\def\F{\mathcal F}
\def\K{\mathcal K}
\def\llen{{\ell\ell}}
\def\O{\mathcal O}
\def\p{{\mathfrak p}}
\let\phi\varphi
\def\Mod{{\bf mod}}
\def\PMod{{\bf proj}}
\def\Pair{\genfrac{}{}{0pt}1}
\def\rad{\operatorname{rad}}
\def\Rad{\operatorname{Rad}}
\def\res{\operatorname{res}}
\def\Ext{\operatorname{Ext}}
\def\sl{\widehat{\mathfrak{sl}}}
\def\soc{\operatorname{Soc}}

\let\part\vdash

\def\iarrow{\stackrel i\longrightarrow}

\newdimen\hoogte    \hoogte=10pt    
\newdimen\breedte   \breedte=15pt   
\newdimen\dikte     \dikte=0.7pt    

\newenvironment{Young}{\begingroup
       \def\vr{\vrule height0.8\hoogte width\dikte depth 0.2\hoogte}
       \def\fbox##1{\vbox{\offinterlineskip
                    \hrule height\dikte
                    \hbox to \breedte{\vr\hfill$##1$\hfill\vr}
                    \hrule height\dikte}}
       \vtop\bgroup\offinterlineskip \tabskip=-\dikte \lineskip=-\dikte
            \halign\bgroup &\fbox{##\unskip}\unskip \crcr}
       {\egroup\egroup\endgroup}

\def\bitab(#1|#2){\Biggl(\ \begin{Young}#1\cr\end{Young}\,,\ 
                           \begin{Young}#2\cr\end{Young}\ \Biggr)}

\begin{document}

\begin{abstract}
This paper determines the representation type of the Iwahori-Hecke
algebras of type $B$ when $q\ne\pm1$. In particular, we show that a
single parameter non--semisimple Iwahori--Hecke algebra of type~$B$ has
finite representation type if and only if $q$ is a simple root of the
Poincar\'e polynomial, confirming a conjecture of
Uno's~\cite{Uno:reptype}.
\end{abstract}

\maketitle


\section{Introduction}

In this paper we determine the representation type of the Hecke
algebras of type~$B$. Previously, Uno~\cite{Uno:reptype} determined
the representation type of the (one parameter) Iwahori--Hecke algebras
for the rank~2 Coxeter groups and the Coxeter groups of type~$A$. We
build upon Uno's work to study the Hecke algebras of type~$B$; in
particular, we settle Uno's conjecture in this case. 

Let $R$ be an integral domain and suppose that $q$ and $Q$ are
invertible elements of~$R$.  The {\sf Iwahori--Hecke algebra}
$\H=\H_{q,Q}(B_n)$ of type  $B_n$ is the unital associative
$R$--algebra with generators $T_0,T_1,\dots,T_{n-1}$ and relations
$$\begin{array}{rl}
  (T_0+1)(T_0-Q)=0,&
  (T_i+1)(T_i-q) =0\text{\ for $1\le i\le n-1$,}\\
  T_0T_1T_0T_1=T_1T_0T_1T_0,&\qquad
  T_{i+1}T_iT_{i+1}=T_iT_{i+1}T_i\text{\ for $1\le i\le n-2$,}\\
 \multicolumn2c{T_iT_j=T_jT_i\text{\ for $0\le i<j-1\le n-2$.}}
\end{array}$$

We will determine the representation type of $\H$.

Let $\H_q(A_{n-1})$ be the subalgebra of $\H$ generated by
$T_1,\dots,T_{n-1}$; then $\H_q(A_{n-1})$ is isomorphic to the
Iwahori--Hecke algebra of the symmetric group of degree $n$. 

Let $e\in\{1,2,3,\dots,\infty\}$ be the multiplicative order of $q$ in
$R$.

\begin{Point}
{(Uno)~\cite[Proposition 3.7,Theorem 3.8]{Uno:reptype}}
Suppose that $R$ is a field and that $q\ne1$. 
Then $\H_q(A_{n-1})$ is of
finite representation type if and only if $n<2e$.
\label{Uno}\end{Point}

Note that although Uno stated the theorem only in the case where $R$ is
the field of complex numbers; this is not essential in his proof. We
also remark that K.~Erdmann and D.~K.~Nakano~\cite[Theorem 1.2]{ED}
have determined the representation type of all of the blocks of
$\H_q(A_{n-1})$; so (\ref{Uno}) also follows from their result. 

The following reduction theorem, together with (\ref{Uno}), will allow
us to assume that $-Q$ is a power of $q$.

\begin{Point}{(Dipper--James) \cite[Theorem 4.17]{DJ:B}} 
Suppose that $Q\ne-q^f$ for any $f\in \Z$. Then $\H_{q,Q}(B_n)$ is Morita
equivalent to 
$$\bigoplus_{m=0}^n\H_q(A_{m-1})\otimes\H_q(A_{n-m-1}).$$
\end{Point}

Combining the last two results we obtain.

\begin{Corollary}
Suppose that $R$ is a field and that 
$Q\ne-q^f$ for any $f\in\Z$.
Then $\H_{q,Q}(B_n)$ is of finite representation type if and only if
$n<2e$.  
\label{reduction}\end{Corollary}

It remains to determine the representation type of $\H$ when $Q=-q^f$
for some $f\in\Z$. When $Q=-q^f$ the relation for $T_0$ becomes
$(T_0+1)(T_0+q^f)=0$. If $e$ is finite we may assume that $0\le f<e$.
It is convenient to renormalize $T_0$ as $-T_0$, when 
$0\le f\le\frac e2$, and as $-q^{-f}T_0$, when $\frac e2<f<e$; in this
way, the relation for $T_0$ becomes $(T_0-1)(T_0-q^f)=0$ where $0\le
f\le\frac e2$ whenever $e$ is finite.

Henceforth we assume that $q$ is a primitive $e^{\text{th}}$
root of unity in $R$ and that~$T_0$ satisfies the relation
$(T_0-1)(T_0-q^f)=0$ where $0\le f\le\frac e2$. As the $R$--algebra $\H$
now only depends on $q$ we now write $\H=\H_q(B_n)$, or
$\H_{R,q}(B_n)$ when we wish to emphasize the choice of $R$.

The main result of this paper is the following. We will consider the
cases $q=\pm1$ (that is, $e=1$ and $e=2$) separately in
\cite{AM:finite}.

\begin{Theorem} Suppose that $R$ is a field and 
that $e\ge3$ and $0\le f\le\tfrac e2$. 
Then $\H_q(B_n)$ is of finite representation type if and only if
$$n<\min\{e,2f+4\}.$$
\label{main}\end{Theorem}

Uno~\cite{Uno:reptype} asked whether the representation type of a
non--semisimple single parameter Iwahori--Hecke algebra is finite if
and only if $q$ is a simple root of the Poincar\'e polynomial of the
corresponding finite Coxeter group. With the assumptions currently in
place, the one parameter Hecke algebra of type $B$ corresponds to~$e$
being even and~$f=\tfrac e2-1$; so our result gives an affirmative
answer to Uno's question in type~$B$. In fact, if $e$ is even and $\H$
is not semisimple then Theorem~\ref{main} says that $\H$ is of finite
representation type if and only if $\tfrac e2=f+1\le n<e$; this is if
and only if~$q$ is a simple root of the Poincar\'e polynomial ($q$ is
a root of the factor~$x^e-1=0$). If $\H$ is a non--semisimple one
parameter Hecke algebra of type $B$ with $e$ odd then, by
Corollary~\ref{reduction}, $\H$ is of finite representation type if
and only if $e\le n<2e$; again, this is if and only if $q$ is a simple
root of the Poincar\'e polynomial (this time~$q$ is a root
of~$x^{2e}-1=0$).

The proof of Theorem~\ref{main} will occupy all of this paper. In
sections~2 and~3 we recall the results that we need from the
representation theory of algebras and from the representation theory of
$\H$; section~4 shows that $\H$ has infinite representation type when
$n\ge e$; section~5 shows that $\H$ has infinite representation
type when $n\ge 2f+4$; finally, section~6 shows that $\H$ has finite
representation type in the remaining cases.


\section{Preliminaries on representation type}

An algebra $A$ has {\sf finite representation type} if there are only
a finite number of isomorphism classes of indecomposable $A$--modules;
otherwise, $A$ has {\sf infinite representation type}. This section
summarizes the results that we need on the representation type of
algebras. More details can be found in the books of Auslander, Reiten
and Smal\o~\cite{ARS} and Benson~\cite{Benson:I}. 

Suppose that $K$ is a field. We always assume that $K$ is a
splitting field for~$A$. The following two results are well--known.
Throughout the paper, all modules are right modules. 

{\samepage
\begin{Lemma}
\label{subalgebra}
Let $A$ be a finite dimensional $K$--algebra.
\begin{enumerate}
\item Suppose that $I$ is a two--sided ideal of $A$ such that $A/I$
has infinite representation type. Then $A$ is of infinite
representation type. 
\item Suppose that $B$ is a direct summand of $A$ as a
$(B,B)$--bimodule. Then
\begin{enumerate}
\item If $B$ is of infinite representation type then so is $A$.
\item If $A$ is of   finite representation type then so is $B$.
\end{enumerate}\end{enumerate}\end{Lemma}
}

\begin{Lemma}
\label{end-finite}
Let $A$ be a finite dimensional $K$--algebra and let 
$P_1,\dots P_l$ be the complete set of projective indecomposable
$A$--modules, up to isomorphism. Then
\begin{enumerate}

\item $A$ is Morita equivalent to $\End_A(P_1\oplus\dots\oplus P_l)$.
\item if $\End_A(P_i)$ has infinite representation type for some 
$i$ then $A$ has infinite representation type.
\item for each $i$ the algebra $\End_A(P_i)$  has finite
representation type if and only if $\End_A(P_i)\cong K[x]/\<x^m\>$ for
some integer $m\ge0$ $($which depends on~$i)$.
\end{enumerate}\end{Lemma}

For any $A$--module $M$ let $\Rad M$ be the Jacobson radical 
of $M$. Let $D_1,\dots,D_l$ be a complete set 
of isomorphism classes of simple $A$--modules and let 
$P_1,\dots,P_l$ be the corresponding projective indecomposables.

In the theory of algebras Dynkin diagrams are valued graphs, 
with the underlying graph being the usual Dynkin diagram; see, for
example, \cite[VII.3, p241]{ARS}.  
If $A$ is a symmetric algebra, then the {\sf separation diagram} 
of $A$ is the valued graph with vertices 
$\{1,\dots,l,1',\dots,l'\}$ and edges 
$i\overset{(a,b)}{\rule[0.6ex]{2em}{.1ex}}j'$ where 
$a=[\Rad P_i/\Rad^2P_i:S_j]$ and $b=[\Rad P_j/\Rad^2P_j:S_i]$. 

The following result is fundamental, and may be derived from 
the theory of hereditary algebras. 

\begin{Theorem}[Gabriel]
Suppose that $A$ is an indecomposable algebra such that $\Rad^2A=0$.
Then $A$ is of finite representation type if and only if the
separation diagram of $A$ is a disjoint union of Dynkin diagrams of
finite type as a valued graph.
\label{Gabriel}\end{Theorem}

The {\sf Auslander--Reiten quiver} of $A$ is the directed graph with
vertices the indecomposable $A$--modules and edges the irreducible
morphisms between the indecomposables (a map $\phi\map MN$ is
irreducible if $\phi$ has no left or right inverse and whenever $\phi$
factorizes as $\phi=\theta\psi$ then either $\theta$ has a right
inverse or $\psi$ has a left inverse).

\begin{Theorem}
[Auslander]
Let $A$ be an indecomposable algebra and suppose that the
Auslander--Reiten quiver of $A$ has a connected component which has a
finite number of vertices. Then $A$ is of finite representation type.
\label{Auslander}\end{Theorem}

Uno used Auslander--Reiten sequences and Theorem~\ref{Auslander} to
prove the following.

\begin{Theorem}
[\protect{Uno~\cite[Theorem~3.6]{Uno:reptype}}]
Suppose that $A$ is a symmetric indecomposable algebra and that the
decomposition matrix of $A$ can be written in the form
$$\left(\begin{array}{*5c} 1 &   0   &\cdots & 0 \\
                 1 &   1   &\cdots & 0 \\
                 \vdots &\ddots &\ddots & \vdots \\
                 0 &\cdots &   1   & 1 \\
                 0 &\cdots &   0   & 1 
        \end{array}\right).$$
Then $A$ is of finite representation type.
\label{weight one}\end{Theorem}

It turns out that in the cases where $\H$ has finite representation
type all of the non--semisimple blocks of $\H$ satisfy the hypotheses
of Theorem~\ref{weight one}; hence they are of finite type. Hence,  in
principle, we can compute all of the indecomposable modules using
Auslander--Reiten sequences when $\H$ has finite type.

We remark that Uno's paper does not actually contain the statement
above. However, the result can be extracted from his paper because the
assumptions of Theorem~\ref{weight one} appear as
\cite[Theorem~3.4]{Uno:reptype} and these are all that are used in the
proof of his Theorem~3.6.


\section{Results from the representation theory of~$\H_q(B_n)$}

We now turn to the representation theory of $\H_q(B_n)$. 

Let $*$ be the anti--involution of $\H_q(B_n)$ determined by
$T_i^*=T_i$ for $0\le i<n$. If $M$ is a right $\H_q(B_n)$-module then
$\Hom_K(M,K)$ is naturally a left $\H_q(B_n)$-module and it becomes a
right module by twisting the $\H_q(B_n)$--action by the 
anti--involution~$*$. We call this the {\sf dual} of $M$; $M$ is 
{\sf self--dual} if it is isomorphic to its dual.

It is well--known that $\H_q(B_n)$ is a symmetric algebra; hence we
have the following.

\begin{Lemma}
The algebra $\H_q(B_n)$ is a symmetric algebra. In particular,
if $D$ is a simple $\H_q(B_n)$--module and $P$ is its projective cover
then $P/\Rad P\cong D$, $\soc P\cong D$ and $P$ is self-dual. 
\label{middle}\end{Lemma}

As $P$ is self-dual, the dual of the radical series of $P$ is the
socle series of $P$. We remark that this does not mean that the
radical series must be symmetric with respect to its middle layer. 

Applying Lemma~\ref{subalgebra}(ii)(a) to the inclusion
$\H_q(B_m)\incl\H_q(B_n)$, for $m\le n$, yields the following.

\begin{Corollary}
Suppose that $m\le n$ and that $\H_q(B_m)$ is of infinite
representation type. Then $\H_q(B_n)$ is of infinite representation
type.
\label{small}\end{Corollary}

Recall that a {\sf partition} of $n$ is an non--increasing sequence
$\sigma=(\sigma_1\ge\sigma_2\ge\dots)$ of non--negative integers
such that $|\sigma|=n$ where $|\sigma|=\sum_i\sigma_i$.  A 
{\sf bipartition} of~$n$ is an ordered pair
$\lambda=(\lambda^{(1)},\lambda^{(2)})$ of partitions $\lambda^{(1)}$
and $\lambda^{(2)}$ such that $|\lambda^{(1)}|+|\lambda^{(2)}|=n$; we
write $\lambda\part n$ and $|\lambda|=n$. The set of bipartitions is
naturally a poset with partial order $\gedom$ where $\lambda\gedom\mu$
if for all $k\ge1$ and $l\ge1$
$$\sum_{i=1}^k\lambda^{(1)}_i\ge\sum_{i=1}^k\mu^{(1)}_i\And
|\lambda^{(1)}|+\sum_{j=1}^l\lambda^{(2)}_j
    \ge|\mu^{(1)}|+\sum_{j=1}^l\mu^{(2)}_j.$$
If $\lambda\gedom\mu$ we say that $\lambda$ {\sf dominates} $\mu$.
If $\lambda\gedom\mu$ and $\lambda\ne\mu$ we write $\lambda\gdom\mu$.

Let $\A=\Z[t,t^{-1}]$ where $t$ is an indeterminate. Then Dipper,
James and Murphy have shown that there exist a family
$\set{S^\lambda_\A|\lambda\part n}$ of free $\A$--modules which are
equipped with operators $T_0,T_1,\dots,T_{n-1}\in\End_\A(S^\lambda)$
which satisfy the defining relations of $\H_t(B_n)$. Consequently,
$S^\lambda=S^\lambda_\A\otimes_\A K$ is a $\H_q(B_n)$--module,
where $\H_q(B_n)$ is the Hecke algebra defined over the field $K$ with
$q\in K$ and we consider~$K$ as an $\A$--module by letting $t$ act on
$K$ as multiplication by $q$.

The module $S^\lambda$ is a {\sf Specht module} of $\H_q(B_n)$. It comes
equipped with a symmetric bilinear form $\<\ ,\ \>$ such that $\<u
T_i,v\>=\<u,vT_i\>$ for $0\le i<n$ and all $u,v\in S^\lambda$.
Consequently, the module 
$$\rad S^\lambda=\set{u\in S^\lambda|\<u,v\>=0\ForAll v\in S^\lambda}$$
is an $\H$--submodule of $S^\lambda.$ Set $D^\lambda=S^\lambda/\rad
S^\lambda$. 

The modules $S^\lambda$ enjoy the following properties. 
Theorem \ref{DJM}(i) is proved by modifying the proof of 
\cite[Theorem 6.1]{DJM}. The others are 
stated in \cite[Theorem 6.5]{DJM} and \cite[Theorem 6.6]{DJM}. 

\begin{Theorem}
[Dipper--James--Murphy~\cite{DJM}]\label{DJM}
Suppose that $K$ is a field.
\begin{enumerate}
\item Any $\H_q(B_n)$--submodule of $S^\lambda$ contains 
$\rad S^\lambda$ or is contained in $\rad S^\lambda$. In particular, 
the module $D^\lambda$ is either $0$ or an absolutely irreducible
self--dual $\H_q(B_n)$--module.
\item $\set{D^\mu|D^\mu\ne0}$ is a complete set of pairwise
non--isomorphic irreducible $\H_q(B_n)$--modules.
\item If $D^\mu\ne0$ then the decomposition multiplicity
$[S^\lambda:D^\mu]\ne0$ only if $\lambda\gedom\mu$; further, 
$[S^\mu:D^\mu]=1.$
\end{enumerate}\end{Theorem}

In particular, if $D^\mu\ne0$ then $D^\mu$ is the unique head of
$S^\mu$ and $\Rad S^\mu=\rad S^\mu$; consequently, $S^\mu$ is
indecomposable and if $S^\mu\ne D^\mu$ then $S^\mu$ has Loewy length
at least $2$. If $D^\mu\ne0$ let $P^\mu$ be the corresponding
principal indecomposable $\H$--module; in other words, $P^\mu$ is the
projective cover of~$D^\mu$. Let $d_{\lambda\mu}=[S^\lambda:D^\mu]$ be
the multiplicity of~$D^\mu$ as a composition factor of~$S^\lambda$.

It is implicit in the work of Dipper, James and Murphy that
$\H_q(B_n)$ is a cellular algebra in the sense of Graham and
Lehrer~\cite{GL} (compare \cite{DJM:cyc}). Consequently, the theory of
cellular algebras gives us the following result.

\begin{Corollary}
Let $P$ be a projective $\H$--module.
\begin{enumerate}
\item Then $P$ has a Specht filtration; thus, there exist
bipartitions $\nu_1,\dots,\nu_k$ and a filtration 
$P=P^k>P^{k-1}>\dots>P^1>0$ such that 
$P^i/P^{i-1}\cong S^{\nu_i}$, for $1<i\le k$, and $i<j$ 
whenever $\nu_i\gdom\nu_j$.
\item Suppose that $P=P^\mu$ for some bipartition $\mu$ with
$D^\mu\ne0$. Then the Specht filtration can be chosen so that
$d_{\lambda\mu}=\#\set{1\le i\le k|\nu_i=\lambda}$. In particular, if
$\lambda$ is maximal in the dominance ordering such that
$d_{\lambda\mu}\ne0$ then $P^\mu$ has a submodule isomorphic to
$S^{\lambda}$. 
\end{enumerate}
\label{specht filtration}
\end{Corollary}

\begin{proof}
These results are implicit in the work of
Graham--Lehrer~\cite{GL} (and slightly more explicit in
\cite[Lemma~2.19]{M:ULect}). The existence of the filtration is
exactly \cite[Lemma~2.9(ii)]{GL}; that we can order the bipartitions
$\nu_i$ by dominance follows from the choice of
$\Phi_0\subset\Phi_1\subset\dots\subset\Phi_d$ as a maximal chain of
ideals in the proof of this result. Part~(ii) follows by combining 
Lemma~2.10(i) and Theorem~3.7(ii) of \cite{GL}.
\end{proof}

We remark that some care must be taken when working with Specht
filtrations because Specht modules indexed by different bipartitions
can be isomorphic when~$\H$ is not semisimple. This technicality can
be avoided by working with a modular system and lifting the projective
module to the discrete valuation ring where the Specht filtration is
unambiguously defined.

In principle, Theorem~\ref{DJM}(ii) produces all of the irreducible
$\H$--modules; however, determining when $D^\lambda$ is non--zero is
still a difficult problem. The non--zero $D^\lambda$ have now been
classified by the first author; to describe this result we need some
more nomenclature.

The {\sf diagram} $[\lambda]$ of a bipartition
$\lambda=(\lambda^{(1)},\lambda^{(2)})$ is the set of {\sf nodes}
$$[\lambda]
    =\set{(i,j,k)|1\le j\le\lambda^{(k)}_i\For i\ge1\And k=1,2},$$ 
which we will think of as being an ordered pair of arrays of boxes in
the plane. 

Given two nodes $x=(i,j,k)$ and $y=(i',j',k')$ we say that $y$ is
{\sf below} $x$ if either $k=k'$ and $i<i'$, or $k<k'$. Further,
$x\in[\lambda]$ is {\sf removable} if $[\lambda]\setminus\{x\}$ is the
diagram of a bipartition; similarly, $y\notin[\lambda]$ is {\sf addable}
if $[\lambda]\cup\{y\}$ is the diagram of a bipartition.  The
{\sf content} of $x$ is $c(x)=j-i+(k-1)f$ and the
{\sf residue} of  $x$ is $\res(x)=c(x)\pmod e$ ---
recall that $q$ is an $e^{\text{th}}$ root of unity and $Q=-q^f$. If
$r=\res(x)$ we call $x$ an $r$--node.

An $r$--node $x$ is {\sf normal} if whenever $y$ is a removable
$r$--node below $x$ then there are more removable $r$--nodes between
$x$ and $y$ than there are addable $r$--nodes, and there are at least
as many removable $r$--nodes below $x$ as addable $r$--nodes below $x$.
In addition, $x$ is {\sf good} if there are no normal $r$--nodes above
$x$. Here, $0\le r<e$.

Finally, a bipartition $\mu$ is {\sf Kleshchev} if either
$\mu=\((0),(0)\)$ or $\mu$ contains a good node $x$ such that
$[\mu]\setminus\{x\}$ is the diagram of a Kleshchev bipartition.

\begin{Theorem}
[Ariki~\cite{A:class}]
Suppose that $\mu$ is a bipartition of $n$. Then $D^\mu\ne0$ if and
only if $\mu$ is a Kleshchev bipartition.
\end{Theorem}

The proof of this result builds on the next theorem which reveals the
deep connections between the representation theory of
$\H_q(B_n)$ and the representation theory of the Kac--Moody algebra
$U(\sl_e)$ of type $A^{(1)}_{e-1}$ in characteristic zero. Let
$\Lambda_0,\dots,\Lambda_{e-1}$ be the fundamental weights of
$U(\sl_e)$ and for each dominant weight $\Lambda$ let $L(\Lambda)$ be
the corresponding integral highest weight module.

Let $\H_q(B_n)$--\Mod\ be the category of finite dimensional right
$\H_q(B_n)$--modules and $\H_q(B_n)$--\PMod\ be the category of finite
dimensional projective $\H_q(B_n)$--modules. Finally, let
$\K_0(\mathcal C)$ be the Grothendieck group of the category 
$\mathcal C$.

\begin{Theorem}
[Ariki~\cite{Ariki:can}]
For $i=0,1,\dots,e-1$ there exist exact functors
$$e_i,f_i\map{\H_q(B_n)\text{--\Mod}}\H_q(B_{n\pm1})\text{--\Mod}$$
such that the operators induced by these and suitably defined operators 
$h_i$ for $i=0,1,\dots,e-1$ and $d$ give
$\K_0=\bigoplus_{n\ge0}\K_0\(\H_q(B_n)\text{--\PMod}\)\otimes_\Z\Q$
the structure of a $U(\sl_e)$--module. Furthermore, 
$\K_0 \cong L(\Lambda_0+\Lambda_f)$ as an $U(\sl_e)$--module and if $K$
is a field of characteristic zero then the principal indecomposable
$\H_q(B_n)$--modules correspond to elements of the Lusztig--Kashiwara
canonical basis of $L(\Lambda_0+\Lambda_f)$ under this isomorphism.
\label{Canonical}\end{Theorem}

A {\sf modular system with parameters} is a modular system
$(K,\O,k)$ such that~$K$ is a field of characteristic zero and $\O$ is a
discrete valuation ring with residue field $k$, together with
parameters $t\in\O\incl K$ and $q\in k$ such that $t$ and $q$ have
the same multiplicative order in $K$ and $k$ respectively and~$t$ maps
to $q$ under the canonical map $\O\to k$. Let $\H_k=\H_{k,q}(B_n)$ and
$\H_R=\H_{R,t}(B_n)$, for $R\in\{K,\O\}$, be the corresponding  Hecke
algebras. Then $\H_R\cong\H_\O\otimes_\O R$ for $R\in\{K,k\}$.

Let $\Lambda=\Lambda_0+\Lambda_f$ and fix a highest weight vector
$v_\Lambda$ in $L(\Lambda)$.

If $\mu$ is a Kleshchev bipartition then we write $P^\mu_R$ for the
principal indecomposable $\H_R$--module for $R\in\{K,\O,k\}$.  So if
$\O$ is a complete discrete valuation ring then
$P^\mu_k=P_\O\otimes_\O k$ and $P^\mu_K$ is a direct summand of
$P_\O\otimes_\O K$. By abuse of notation,  we also let $[P^\mu_K]$
denote both the equivalence class of $P^\mu_K$ in $\K_0$ and the
corresponding canonical basis element of $L(\Lambda)$.

\begin{Corollary}
Let $(K,\O,k)$ be a modular system with parameters. Suppose
that~$\mu$ is a Kleshchev bipartition such that 
$[P^\mu_K]=f_{i_1}^{(m_1)}\dots f_{i_l}^{(m_l)}v_\Lambda$
for some $m_1,\dots,m_l$ and $i_1,\dots,i_l$. Then the decomposition
map sends $[P^\mu_K]$ to $[P^\mu_k]$.
\label{canonical}
\end{Corollary}

\begin{proof}
Let $\O\subset \tilde{\O}$ be an embedding into a complete discrete
valuation ring $\tilde O$ where~$\tilde O$ has residue field $k$ and
$K\subset \tilde K$  is a field extension such that
$\tilde\O\subset\tilde K$.  Since $D^\mu$ is absolutely irreducible,
$P_{\tilde K}^\mu=P_K^\mu\otimes \tilde K$.  Hence it is enough to
prove the statement under the assumption that $\O$ is a complete
discrete valuation ring.

Set $N=m_1!\dots m_l!$ and let 
$M_\O=f_{i_1}^{m_1}\dots f_{i_l}^{m_l}1_{\H_\O}$ where $1_{\H_\O}$ is
the trivial $\H_{\O,t}(B_0)$--module. Then $M_\O$ is a projective
$\H_\O$--module by the definition of the~$f_i$. 
Therefore, by Corollary~\ref{specht filtration} applied to $\H_\O$,
the module $M_\O$ has a Specht
filtration. Let $M_K=M_\O\otimes K$; then, by assumption, 
$[M_K]=N[P^\mu_K]$; therefore, 
$[M_K]=N[S^\mu_K]+\sum_{\lambda\gdom\mu}N d_{\lambda\mu}[S^\lambda_K]$, 
where $d_{\lambda\mu}=[S^\lambda_K:D^\mu_K]$. Consequently, there is
a surjective homomorphism $M_\O\to{S^\mu_\O}^{\oplus N}$;
tensoring with $k$ gives a surjective homomorphism
$M_k\overset{\phi}\to{D^\mu_k}^{\oplus N}$ (as $D^\mu$ is the head
of $S^\mu$). Since~$M_k$ is a
projective $\H_k$--module there exists a map $\psi$ which makes the
following diagram commute.
$$\psmatrix[colsep=1cm,rowsep=1cm]
                    & M_k\\
  {P^\mu_k}^{\oplus N}& {D^\mu_k}^{\oplus N} & 0
  \everypsbox{\scriptstyle}
  \psset{arrows=->,nodesep=3pt}
  \ncline{1,2}{2,2}\trput{\phi}
  \ncline[linestyle=dashed]{1,2}{2,1}\taput{\psi}
  \ncline{2,1}{2,2}
  \ncline{2,2}{2,3}
\endpsmatrix$$
Now, $\phi$ is surjective so $\psi$ must also be surjective.  On the
other hand, $P^\mu_k$ is a projective $\H_k$--module so we also have a
commutative diagram
$$\psmatrix[colsep=1cm,rowsep=1cm]
       & {P^\mu_k}^{\oplus N}\\
   M_k & {P^\mu_k}^{\oplus N}& 0
  \everypsbox{\scriptstyle}
  \psset{arrows=->,nodesep=3pt}
  \ncline[doubleline=true,doublesep=.8mm]{-}{1,2}{2,2}
  \ncline[linestyle=dashed]{1,2}{2,1}
  \ncline{2,1}{2,2}\tbput{\psi}
  \ncline{2,2}{2,3}
\endpsmatrix$$
Hence, ${P^\mu_k}^{\oplus N}$ is a direct summand of $M_k$.
Note that 
\begin{align*}
\dim_k{P_k^\mu}^{\oplus N}
    &=\dim_K(P_\O^\mu\otimes K)^{\oplus N}
     \ge\dim_K{P_K^\mu}^{\oplus N}\\
    &=\dim_K M_K=\dim_kM_k 
     \ge\dim_k{P_k^\mu}^{\oplus N};
\end{align*}
thus, $M_k={P^\mu_k}^{\oplus N}$.  However, 
$M_K={P^\mu_K}^{\oplus N}$; so we have shown that the modular 
reduction of ${P^\mu_K}^{\oplus N}$ is~${P^\mu_k}^{\oplus N}$, as 
required.
\end{proof}

In order to apply the last two results we need to set up the machinery
for computing the canonical basis elements $[P^\lambda_K]$.
Let $v$ be an indeterminant over $\Z$ and let $\A=\Z[v,v^{-1}]$.
The {\sf Fock space} is the infinite dimensional $\A$--module
$$\F_\A=\bigoplus_{n\ge0}\bigoplus_{\lambda\part n}\A\lambda.$$

Let $U_\A(\sl_e)$ be Lusztig's $\A$--form of the quantum group of
$U(\sl_e)$. Then there are $v$--analogues $E_i$ and $F_i$ of the
operators $e_i$ and $f_i$ which act on $\F_\A$ and give it the
structure of a $U_\A(\sl_e)$--module; an explicit description of $E_i$
and $F_i$ is given below. The $U_\A(\sl_e)$--submodule of~$\F_\A$
generated by the bipartition $\((0),(0)\)$ is isomorphic to
$L_\A(\Lambda)$, the $\A$--form of $L(\Lambda)$. (Recall that
$\Lambda=\Lambda_0+\Lambda_f$.)

Identifying $L_\A(\Lambda)$ with $U_\A(\sl_e)\cdot\((0),(0)\)$,
Theorem~\ref{Canonical} can be reinterpreted as saying that if $\mu$
is a Kleshchev bipartition then there exist polynomials
$d_{\lambda\mu}(v)$ such that 
$$[P^\mu_K]=\mu+\sum_{\lambda\part n}d_{\lambda\mu}(v)\lambda,
\quad\Where d_{\lambda\mu}(v)\in v\Z[v],\Number{canonical basis}$$ 
and $d_{\lambda\mu}=d_{\lambda\mu}(1)$. Uglov~\cite{Uglov} has given
an explicit algorithm for computing the canonical basis elements
$[P^\lambda_K]$. Uglov actually works with a different
Fock space; however, we can compute the canonical basis of inside our
Fock space by modifying his algorithm. For the applications we have in
mind it is enough to know that if 
$F_{i_1}^{(m_1)}\dots F_{i_l}^{(m_l)}\cdot\((0),(0)\)$ can be written
in the form of the right hand side of (\ref{canonical basis}) then it
is an element of the canonical basis of~$L_\A(\Lambda)$; hence, in
such situations we may apply Corollary~\ref{canonical}. We now recall
from \cite{AM:simples} how $U_\A(\sl_e)$ acts on $\F_\A$.

Suppose that $\lambda$ is a bipartition of $n-1$ and that $\mu$ is a
bipartition of $n$. We write $\lambda\iarrow\mu$ if
$[\mu]\setminus[\lambda]=\{x\}$ and $\res(x)=q^i$. For $0\le i<e$
let
\begin{align*}
N^r_i(\nu,\lambda)&=\#\set{\nu\iarrow\alpha|\alpha\gdom\lambda}
         -\#\set{\beta\iarrow\lambda|\nu\gdom\beta},\\
N^l_i(\lambda,\mu)&=\#\set{\lambda\iarrow\alpha|\mu\gdom\alpha}
         -\#\set{\beta\iarrow\mu|\beta\gdom\lambda}.
\end{align*}
By \cite[Prop.~2.6]{AM:simples} the action of $E_i$ and $F_i$ on
$\F_\A$ is determined by
$$E_i\lambda=\Sum_{\nu\iarrow\lambda} v^{-N^r_i(\nu,\lambda)}\nu,
\qquad\And\qquad
F_i\lambda=\Sum_{\lambda\iarrow\mu} v^{N^l_i(\lambda,\mu)}\mu.$$
We will only need the formula giving the action of $F_i$. There are
similar formulae for the action of the remaining generators of
$U_\A(\sl_e)$.

The final tool that we shall need is the analogue of the Jantzen sum
formula for Hecke algebras of type $B$ over an arbitrary field $K$.
The setup is a little technical; we include it for completeness. Let
$\p$ be the maximal ideal of $K[t]$ generated by $t-q$ and let
$\O=K[t]_{\p}$, where $t$ is an indeterminate over $K$; then
$K\cong\O/\p$ (so $\O$ is a localized ring and, in particular, a
discrete valuation ring). Let $\H_\O$ be the Hecke algebra over $\O$
with parameters $t$ and $-q^f(t-q+1)^n$; then $\H_\O\otimes K(t)$ is
semisimple and $\H_K$ is the reduction of $\H_\O$ modulo $\p$. Let
$\nu_\p$ be the $\p$--adic valuation on $\O$.

Previously we defined the residue $\res(x)$ of a node $x=(i,j,k)$.
Define the {\sf $\O$--residue} of $x=(i,j,k)$ to be
$\res_\O(x)=t^{j-i}q^{(k-1)f}(t-q+1)^{(k-1)n}$. The relationship
between these two definitions is that $\res_\O(x)\otimes
1_K=q^{\res(x)}$.

Let $\lambda$ be a bipartition and for each node
$x=(i,j,k)\in[\lambda]$ let $r_x$ be the corresponding rim hook (so
$r_x$ is a rim  hook in $[\lambda^{(k)}]$); the point is that
$[\lambda]\setminus r_x$ is the diagram of a bipartition. Let
$\llen(r_x)$ be the leg length of $r_x$ and define
$\res_\O(r_x)=\res_\O(f_x)$ where $f_x$ is the foot node of $r_x$. The
definitions of these terms can be found, for example, in
\cite{M:ULect}.

Suppose that $\lambda$ and $\mu$ are bipartitions of $n$. If
$\lambda\not\gdom\mu$ let $g_{\lambda\mu}=1$; otherwise set
$$g_{\lambda\mu}=\prod_{x\in[\lambda]} 
      \prod_{\stackrel{y\in[\mu]}
            {[\mu]\setminus r_y=[\lambda]\setminus r_x}}
             \(\res_\O(r_x)-\res_\O(r_y)\)^{\varepsilon_{xy}},$$
where $\varepsilon_{xy}=(-1)^{\llen(r_x)+\llen(r_y)}$. The
$g_{\lambda\mu}$ are not as complicated as their definition suggests;
they have a nice combinatorial interpretation, see
\cite[Example~3.39]{JM:cyc-Schaper}.

Finally, let $S^\lambda_\O$ and $S^\lambda_K$ be the Specht modules
for $\H_\O$ and $\H_K$ respectively. For each $i\ge0$ define 
$S^\lambda_\O(i)
  =\set{u\in S^\lambda_\O|\<u,v\>\in\p^i\ForAll v\in S^\lambda_\O}$.
Then the {\sf Jantzen filtration} of $S^\lambda_K$ is the filtration
$$S^\lambda_K=S^\lambda_K(0)\ge S^\lambda_K(1)
              \ge S^\lambda_K(2)\ge\dots$$
where 
$S^\lambda_K(i)
      =\(S^\lambda_\O(i)+\p S^\lambda_\O\)/\p S^\lambda_\O$. 
In particular, note that $\rad S^\lambda_K=S^\lambda_K(1)$.

We can now state the analogue of Jantzen's sum formula for $\H_K$.

\begin{Theorem}
[\protect{James--Mathas~\cite[Theorem~4.6]{JM:cyc-Schaper}}]
Let $\lambda$ be a bipartition of $n$. Then
$$\sum_{i>0}[S^\lambda_K(i)]
   =\sum_{\mu:\lambda\gdom\mu}\nu_\p(g_{\lambda\mu})[S^\mu_K].$$
in the Grothendieck group $\K_0(\H_K\text{--\Mod})$ of $\H_K$.
\label{sum formula}\end{Theorem}

In general, if $\lambda\gdom\mu$ then $\nu_\p(g_{\lambda\mu})$ is
non--zero only if it is possible to remove a rim hook $r_x$
from $\lambda$ and reattach it to $\mu$ without changing the residue
$\res(r_x)$ of the foot node. In fact, we will only apply this result
when $n<e$; in this situation we have
$\nu_\p(g_{\lambda\mu})\in\{0,1\}$ so the technicalities above can be
ignored.

All of the composition factors of a Specht module belong to the same
block (for example, because $\H$ is cellular); we abuse notation and
say that $\lambda$ is contained in the block~$B$ if $S^\lambda$ is
contained in $B$. Say that two bipartitions $\lambda$ and $\mu$ are
{\sf linked by hooks} if there is a sequence of bipartitions
$\lambda=\nu_1,\dots,\nu_l=\mu$ such that, for each~$i$,
$[\nu_{i+1}]\setminus r_{y_i}=[\nu_i]\setminus r_{x_i}$ and
$\res(r_{x_i})=\res(r_{y_i})$ for some nodes $x_i\in[\nu_{i+1}]$
and $y_i\in[\nu_i]$.

\begin{Proposition}
Suppose that $S^\lambda$ and $S^\mu$ are in the same block. Then
$\lambda$ and $\mu$ are linked by hooks. 
\label{cor of sum formula}\end{Proposition}

\begin{proof}
By definition, $S^\lambda$ and $S^\mu$ are in the same block if and
only if there exists a sequence of bipartitions
$\mu=\nu_1,\dots,\nu_l=\lambda$ such that $S^{\nu_i}$ and
$S^{\nu_{i+1}}$ have a common composition factor. Thus, it is enough
to prove that if $D^\mu\ne0$ appears in $S^\lambda$ then~$\lambda$ and
$\mu$ are linked by hooks. If $\lambda\ne\mu$ then
$\lambda\gdom\mu$ by Theorem~\ref{DJM}(iii). The sum formula
implies that $D^\mu$ appears in $S^{\nu}$ for some $\nu$ such that
$\lambda\gdom\nu$ and $\lambda$ and $\nu$ are linked by hooks. By
induction on dominance $\nu$ and $\mu$ are linked by hooks so we are
done. 
\end{proof}

If $\lambda$ is a bipartition let $\res(\lambda)$ be the {\it
multiset} $\set{\res(x)|x\in[\lambda]}$. Then as a corollary of the
Proposition we have the following (there is an easier proof).

\begin{Corollary}
[\protect{Dipper--James--Murphy~\cite[Corollary 8.7]{DJM}}]
Suppose that $S^\lambda$ and~$S^\mu$ are in the same block. Then
$\res(\lambda)=\res(\mu)$ $($as multisets$)$.
\label{blocks}\end{Corollary}

By the Corollary we can define the {\sf residue } of a block~$B$ to
be the multiset $\res(B)=\res(\lambda)$ where $\lambda$ is any
bipartition  contained in~$B$. 

In fact, Grojnowski~\cite{Groj:AKblocks} has recently shown that the
converse of Corollary~\ref{blocks} is true; so, two Specht modules
$S^\lambda$ and $S^\mu$ belong to the same block if and only if
$\res(\lambda)=\res(\mu)$. We will not need this stronger result. 


\section{The representation type when $n=e$}

In this section we will prove the following result. Recall that we are
assuming that $e\ge3$.

\begin{Theorem}
Suppose $n\ge e$. Then $\H$ has infinite
representation type.
\label{n=e}
\end{Theorem}

\begin{proof}
By Corollary~\ref{small} we may assume that $n=e$. Further, by
Lemma~\ref{subalgebra} it is enough to show that one block of $\H$ has
infinite representation type; we will show that the block $B$ with
residues $\{0,1,\dots,e-1\}$ has infinite representation type. 

There will be several cases to consider. To begin suppose that
$f\ne0$. Because all of the residues in $B$ are distinct a bipartition
$\lambda=(\lambda^{(1)},\lambda^{(2)})$ appears in $B$ only
if~$\lambda^{(1)}$ and $\lambda^{(2)}$ are both hook partitions; that
is, $\lambda=\((a,1^b),(c,1^d)\)$ for some $a,b,c,d$. Define
bipartitions
\begin{xalignat*}{2}
  \lambda_k&=\((0),(k,1^{e-k})\), &&\For 1\le k\le e,\\
  \mu_k&=\((k,1^{e-k}),(0)\),&&\For 1\le k\le e,\\
  \lambda_{k,l}&=\((f-l,1^{e-f-k}),(k,1^l)\),&&
           \For 1\le k\le e-f\And 0\le l<f.
\end{xalignat*}
It is easily checked that all of these bipartitions belong to $B$.
Certainly, the two sets of bipartitions $\{\lambda_k\}$ and
$\{\mu_k\}$ are disjoint; the restrictions on $k$ and $l$ ensure that
$\lambda_{k,l}\ne\lambda_m$ and $\lambda_{k,l}\ne\mu_m$ for any $m$.
Consequently, this is a complete list of the bipartitions which appear
in $B$, with no repeats. 

\begin{Proposition}
\label{f>0}
Suppose that $n=e$ and that $1\le f\le\frac e2$.
\begin{enumerate}
\item The complete set of Kleshchev bipartitions in $B$ is
$$\set{\lambda_k|1\le k<e}
  \cup\set{\lambda_{k,l}|1\le k\le e-f\And 0\le l<f}.$$
\item For $1\le k<e$ we have
$$[P^{\lambda_k}]=[S^{\lambda_k}]+[S^{\lambda_{k+1}}]
    +\begin{cases}%
 [S^{\lambda_{k,f-1}}]+[S^{\lambda_{k+1,f-1}}],& \If k<e-f,\\\relax
 [S^{\lambda_{e-f,f-1}}],&\If k=e-f,\\\relax
 [S^{\lambda_{e-f,e-k}}]+[S^{\lambda_{e-f,e-k-1}}],& \If k>e-f.
\end{cases}$$
\item For $1\le k\le e-f$ and $0\le l<f$ we have
\end{enumerate}
$$[P^{\lambda_{k,l}}]=[S^{\lambda_{k,l}}]+
\begin{cases} [S^{\lambda_{k-1,0}}]+[S^{\mu_{f+k-1}}]+[S^{\mu_{f+k}}],
         &\text{if\ } k\ne1\And l=0,\\\relax
[S^{\mu_f}]+[S^{\mu_{f+1}}],
         &\text{if\ } k=1\And l=0,\\\relax
[S^{\lambda_{k,l-1}}]+[S^{\lambda_{k-1,l}}]+[S^{\lambda_{k-1,l-1}}],
         &\text{if\ } k\ne1\And l\ne0,\\\relax
[S^{\lambda_{1,l-1}}]+[S^{\mu_{f-l}}]+[S^{\mu_{f-l+1}}],&
    \text{if\ } k=1\And l\ne0.
\end{cases}$$
\end{Proposition}

\begin{proof}
It is easy to see that the bipartitions $\lambda_e$ and
$\mu_k$, for $1\le k\le e$, are not Kleshchev. Next suppose that
$1\le k\le e-f$. Then a straightforward computation shows that 
$$\begin{array}{ll}
\multicolumn2l
{F_{f+k}\dots F_{e-1}F_0F_{f+k-1}\dots F_{f+1}F_1\dots F_f\((0),(0)\)}
 \\\qquad
 &=F_{f+k}\dots F_{e-1}F_0\((0),(k,1^{f-1})\)\\
 &=\begin{cases}
  \lambda_k+v\lambda_{k+1} +v\lambda_{k+1,f-1} +v^2\lambda_{k,f-1},
         &\If 1\le k<e-f,\\
  \lambda_{e-f}+v\lambda_{e-f+1}+v^2\lambda_{e-f,f-1},
         &\If k=e-f.
\end{cases}
\end{array}$$
Note that $\((0),(k,1^{f-1})\)$ has two addable $0$--nodes, 
an addable $f$--node and an addable $f+k$--node which is a $0$--node 
if $k=e-f$, and when we add an addable $r$--node, there is no 
removable $r$--node. 

This shows that $\lambda_k$ is Kleshchev for $1\le k\le e-f$; the
formula for $[P^{\lambda_k}]$ now follows from
Corollary~\ref{canonical} (and the remarks after 
(\ref{canonical basis})).

The remaining cases are similar: for $\lambda_k$ with $e-f<k<e$ 
we have $f\ge 2$ and we compute 
$$\begin{array}{ll}
\multicolumn2l{F_{k-e+f}\dots F_{f-1}F_{k-e+f-1}\dots F_1F_0
F_{e-1}\dots F_{f+1}F_f\((0),(0)\)}\\
\qquad&=F_{k-e+f}\dots F_{f-1}F_{k-e+f-1}\dots F_1F_0\((0),(e-f)\)\\
&=F_{k-e+f}\dots F_{f-1}\Big[\((0),(k)\)+v\((k-e+f),(e-f)\)\Big]\\
&=F_{k-e+f}\Big[\((0),(k,1^{e-k-1})\)+v\((k-e+f),(e-f,1^{e-k-1})\)\Big]\\
&=\lambda_k+v\lambda_{k+1}+v\lambda_{e-f,e-k}+v^2\lambda_{e-f,e-k-1};
\end{array}$$
for $\lambda_{k,l}$ with $1\le k\le e-f$ and $l=0$ compute
$$\begin{array}{ll}
\multicolumn2l{F_{f+k-1}\dots F_fF_{f+k}\dots F_{e-1}F_{f-1}
        \dots F_0\((0),(0)\)}\\\qquad
&=F_{f+k-1}\dots F_fF_{f+k}\dots F_{e-1}\((f),(0)\)\\
&=F_{f+k-1}\dots F_f\((f,1^{e-f-k}),(0)\)\\
&=\begin{cases}
  \lambda_{k,0}+v\lambda_{k-1,0} +v\mu_{f+k-1} +v^2\mu_{f+k},
         &\If k\ne1,\\
  \lambda_{1,0}+v\mu_f+v^2\mu_{f+1},
         &\If k=1.
\end{cases}
\end{array}$$
for $\lambda_{k,l}$ with $0<l<f$ and $2\le k\le e-f$ compute
$$\begin{array}{ll}
\multicolumn2l{F_{f-l}\dots F_{f-1}F_{f+k-1}\dots F_f
F_{f+k}\dots F_{e-1}F_{f-l-1}\dots F_0\((0),(0)\)}\\\qquad
&=F_{f-l}\dots F_{f-1}F_{f+k-1}\dots F_f\((f-l,1^{e-f-k}),(0)\)\\
&=F_{f-l}\Big[\((f-l,1^{e-f-k}),(k,1^{l-1})\)\\&
\qquad\qquad+v\((f-l,1^{e-f-k+1}),(k-1,1^{l-1})\)\Big];
\end{array}$$
and, finally, for $\lambda_{1,l}$ with $0<l<f$ compute
$$\begin{array}{ll}
\multicolumn2l{F_{f-l}\dots F_{f-1}F_f
F_{f+1}\dots F_{e-1}F_{f-l-1}\dots F_0\((0),(0)\)}\\\qquad
&=F_{f-l}\dots F_{f-1}F_f\((f-l,1^{e-f-1}),(0)\)\\
&=F_{f-l}\dots F_{f-1}
\Big[\((f-l,1^{e-f-1}),(1)\)+v\((f-l,1^{e-f}),(0)\)\Big]\\
&=F_{f-l}\Big[\((f-l,1^{e-f-1}),(1^l)\)+v\((f-l,1^{e-f+l-1}),(0)\)\Big].
\end{array}$$
In each case, an application of \ref{canonical} now completes the
proof.
\end{proof}

By Lemma~\ref{end-finite}(i) in order to prove Theorem~\ref{n=e} it
is enough to show that $\End_\H(P^{\lambda_1})$ is not isomorphic to
$k[x]/\<x^m\>$ for any $m$. We need to consider several cases.
First we observe that 
$$[P^{\lambda_1}]=[S^{\lambda_1}]+[S^{\lambda_2}]
         +[S^{\lambda_{2,f-1}}]+[S^{\lambda_{1,f-1}}]\Number{PIM}$$
by Proposition~\ref{f>0}(i). We will use this to determine the structure of
$P^{\lambda_1}$.

\Case{$e\ge 5$ and $f\ge 2$}
Then $e-f\ge3$ since $e-f\ge\frac e2=2.5$. By (\ref{canonical basis})
we can use Proposition~\ref{f>0} to compute the decomposition numbers
$d_{\lambda\mu}$ for the four Specht modules appearing in (\ref{PIM});
this gives the following table (omitted entries are zero).
$$\begin{array}{l|*5c}
&D^{\lambda_1}&D^{\lambda_2}&D^{\lambda_{3,f-1}}
                &D^{\lambda_{2,f-1}}&D^{\lambda_{1,f-1}}\\\hline
S^{\lambda_1}&         1 & . & . & . & .\\
S^{\lambda_2}&         1 & 1 & . & . & .\\
S^{\lambda_{2,f-1}}&   1 & 1 & 1 & 1 & .\\
S^{\lambda_{1,f-1}}&   1 & 0 & 0 & 1 & 1\\
\end{array}$$
Consequently, 
$[P^{\lambda_1}]=4[D^{\lambda_1}]+2[D^{\lambda_2}]
         +[D^{\lambda_{1,f-1}}]+2[D^{\lambda_{2,f-1}}]
         +[D^{\lambda_{3,f-1}}]$.
By Corollary~\ref{specht filtration}(ii), $P^{\lambda_1}$ has submodule
isomorphic to $S^{\lambda_{1,f-1}}$. Therefore, $S^{\lambda_{1,f-1}}$
has both a simple head and simple socle; so, looking at the submatrix
of the decomposition matrix above, the Loewy structure of
$S^{\lambda_{1,f-1}}$ is
$$S^{\lambda_{1,f-1}}=\begin{array}{c}D^{\lambda_{1,f-1}}\\
           D^{\lambda_{2,f-1}}\\D^{\lambda_1}  \end{array}.$$ 
We also have 
$$S^{\lambda_2}=\begin{array}{c}D^{\lambda_2}\\
           D^{\lambda_1}  \end{array}.$$ 
Considering the dual of $\Rad S^{\lambda_{1,f-1}}$ and $S^{\lambda_2}$, 
we conclude that $D^{\lambda_2}$ and $D^{\lambda_{2,f-1}}$ appear in 
$\Rad P^{\lambda_1}/\Rad^2P^{\lambda_1}$. If 
$0 \To D^{\lambda_1} \To X \To D^{\lambda_1} \To 0$  is an exact
sequence then, for $1\le i<n$, $T_i-q$ acts invertibly on $X$ so that
$T_i$ acts as~$-1$ on~$X$. Similarly, $T_0$ also acts as
multiplication by a scalar on $X$ since $f\ne0$. Therefore, every such
exact sequence splits and so $\Ext^1(D^{\lambda_1},D^{\lambda_1})=0$;
consequently, $D^{\lambda_1}$ is not a composition factor of~$\Rad
P^{\lambda_1}/\Rad^2P^{\lambda_1}$. 

Again, by Corollary~\ref{specht filtration} $P^{\lambda_1}$ has a
Specht filtration so $\Rad P^{\lambda_1}/S^{\lambda_{1,f-1}}$ has a
Specht filtration whose successive quotients are $S^{\lambda_2}$ and
$S^{\lambda_{2,f-1}}$.  This means that $\Rad
P^{\lambda_1}/\Rad^2P^{\lambda_1}= D^{\lambda_2}\oplus
D^{\lambda_{2,f-1}}$. 

We now use the fact that 
$\Rad S^{\lambda_{1,f-1}}
  =\Pair{D^{\lambda_{2,f-1}}}{D^{\lambda_1}}$ to prove that 
$D^{\lambda_1}$ is contained in 
$\Rad S^{\lambda_{2,f-1}}/\Rad^2 S^{\lambda_{2,f-1}}$. 
As $\Rad S^{\lambda_{1,f-1}}$ has a unique head there is a surjection 
$P^{\lambda_{2,f-1}}\To\Rad S^{\lambda_{1,f-1}}$. Note also that
$P^{\lambda_{2,f-1}}$ has a Specht filtration with successive
quotients $S^{\lambda_{2,f-1}}$ and $S^{\lambda_{2,f-2}}$,
$S^{\lambda_{1,f-1}}$, $S^{\lambda_{1,f-2}}$. 
Since $S^{\lambda_{2,f-2}}$, $S^{\lambda_{1,f-1}}$, $S^{\lambda_{1,f-2}}$ 
must map to $D^{\lambda_1}$, and each of these 
has unique head which 
is not isomorphic to $D^{\lambda_1}$, this surjection induces a map
$S^{\lambda_{2,f-1}}\To\Rad S^{\lambda_{1,f-1}}$. Therefore, 
$\Rad S^{\lambda_{2,f-1}}/\Rad^2 S^{\lambda_{2,f-1}}$ contains
$D^{\lambda_1}$ as a summand. Let $U$ be a module such that
$\Rad^2 S^{\lambda_{2,f-1}}\subseteq U\subseteq\Rad S^{\lambda_{2,f-1}}$
and $\Rad S^{\lambda_{2,f-1}}/U\cong D^{\lambda_1}$ and set
$$V=(\Rad P^{\lambda_1})/(U+S^{\lambda_{1,f-1}}).$$ 
Then there is a short exact sequence
$$0\To\Pair{D^{\lambda_{2,f-1}}}{D^{\lambda_1}}\To V
      \To\Pair{D^{\lambda_2}}{D^{\lambda_1}}\To 0$$
and $V/\Rad V\cong D^{\lambda_2}\oplus D^{\lambda_{2,f-1}}$.
Hence, we also have
$0\To D^{\lambda_1}\To \Rad V\To D^{\lambda_1}\To0$
and consequently $\Rad V=D^{\lambda_1}\oplus D^{\lambda_1}$.
Therefore, $\Rad^2 P^{\lambda_1}/\Rad^3 P^{\lambda_1}$ contains
$D^{\lambda_1}\oplus D^{\lambda_1}$ and it follows that
$P^{\lambda_1}/\Rad^3P^{\lambda_1}$ contains $D^{\lambda_1}\oplus
D^{\lambda_1}$ as an $\H$--submodule, that another $D^{\lambda_1}$
appears as the head of $P^{\lambda_1}/\Rad^3P^{\lambda_1}$, and that
these are the only ways in which $D^{\lambda_1}$ appear as a
composition factor of $P^{\lambda_1}/\Rad^3 P^{\lambda_1}.$ 
Therefore, we have that
$\End_{\H/\Rad^3\H}(P^{\lambda_1}/\Rad^3 P^{\lambda_1})
       \not\cong k[x]/\<x^m\>$. 
Hence, $\End_\H(P^{\lambda_1})$ is not of finite representation type
by Lemma~\ref{end-finite}, so $\H$ has infinite representation type by
Lemma~\ref{subalgebra}(ii).

\Case{$e=4$ and $f=2$} By Proposition~\ref{f>0} and (\ref{PIM}), 
we have the same table as above if we replace $\lambda_{3,f-1}$ by 
$\lambda_3$. Thus 
$$[P^{\lambda_1}]=4[D^{\lambda_1}]+2[D^{\lambda_2}]+[D^{\lambda_3}]
                  +2[D^{\lambda_{2,1}}]+[D^{\lambda_{1,1}}].$$
The argument used in case~1 shows that $P^{\lambda_1}/\Rad^3P^{\lambda_1}$ 
contains $D^{\lambda_1}\oplus D^{\lambda_1}$ as an $\H$--submodule. 
Hence, $\End_\H(P^{\lambda_1})$ is again of infinite type.

\Case{$f=1$} Again using Proposition~\ref{f>0} and (\ref{PIM}) we find that
$$[P^{\lambda_1}]=4[D^{\lambda_1}]+2[D^{\lambda_2}]+[D^{\lambda_{1,0}}]
             +2[D^{\lambda_{2,0}}]+\delta_{e\ge4}[D^{\lambda_{3,0}}],$$
where $\delta_{e\ge4}=1$ if $e\ge4$ and $\delta_{e\ge4}=0$ otherwise.
Almost the same argument as before again shows that 
$P^{\lambda_1}/\Rad^3P^{\lambda_1}$ contains 
$D^{\lambda_1}\oplus D^{\lambda_1}$ as an $\H$--submodule; however,
to show that $\Rad S^{\lambda_{2,0}}/\Rad^2 S^{\lambda_{2,0}}$
contains $D^{\lambda_1}$ 
we need to argue in a different way: since 
$[P^{\lambda_{2,0}}]=[S^{\lambda_{2,0}}]+[S^{\lambda_{1,0}}]
+[S^{\mu_2}]+[S^{\mu_3}]$ and 
$[S^{\mu_2}]=[D^{\lambda_{2,0}}]+[D^{\lambda_{1,0}}]$, 
$[S^{\mu_3}]=[D^{\lambda_{3,0}}]+[D^{\lambda_{2,0}}]$ do not 
contain $D^{\lambda_1}$, the surjection 
$P^{\lambda_{2,0}}\To\Rad S^{\lambda_{1,0}}$ induces 
a surjection $S^{\lambda_{2,0}}\To\Rad S^{\lambda_{1,0}}$. 
Thus, once again $\H$ is of infinite representation type
by Lemma~\ref{end-finite} and Lemma~\ref{subalgebra}.

\Case{$f=0$} It remains to consider the case $f=0$. This case is
somewhat degenerate as $1$ is the only possible eigenvalue for the
action of $T_0$ upon a representation. Because $f=0$ there are no
bipartitions of the form $\lambda_{k,l}$ in $B$ so Proposition~\ref{f>0} cannot
hold in this case. A similar argument shows that we have the following
simpler statement when $f=0$.

\begin{Proposition}
\label{f=0} Suppose that $f=0$. 
\begin{enumerate}
\item $\set{\lambda_k|1\le k<e}$ is the complete set of Kleshchev
bipartitions in $B$.
\item For $1\le k<e$ we have
$[P^{\lambda_k}]=[S^{\lambda_k}]+[S^{\lambda_{k+1}}]
                      +[S^{\mu_k}]+[S^{\mu_{k+1}}].$
\end{enumerate}\end{Proposition}

To apply this result we do not have to work as hard as in the previous
cases. First observe that for $1\le k<e$ the module $D^{\lambda_k}$
can be constructed by letting~$T_0$ act as $1$ on the simple
$\H(A_e)$--module $D^{(k,1^{e-k})}$ (by induction both modules have
the same dimension; namely, $\binom{e-1}{k-1}$). Let $M^{\lambda_k}$
be the $\H$--module obtained by letting~$T_0$ act as
$\big(\begin{smallmatrix}1&1\\0&1\end{smallmatrix}\big)$ on the
$\H(A_e)$--module $D^{(k,1^{e-k})}\oplus D^{(k,1^{e-k})}$. Then
$M^{\lambda_k}$ cannot be semisimple as an $\H$--module because $T_0$
does not act as a scalar. On the other hand, the socle
of~$M^{\lambda_k}$ is simple (being isomorphic to
$D^{\lambda_k}$); hence,~$M^{\lambda_k}$ is indecomposable.  This
implies that $[\Rad P^{\lambda_k}/\Rad^2
P^{\lambda_k}:D^{\lambda_k}]\ne0$.

Now $[S^{\lambda_k}]=[D^{\lambda_k}]+[D^{\lambda_{k-1}}]$ for $1\le
k<e$ by Proposition~\ref{f=0}. Hence, $S^{\lambda_k}$ is
indecomposable and
$[\Rad P^{\lambda_k}/\Rad^2 P^{\lambda_k}:D^{\lambda_{k-1}}]\ne0$
for $2\le k<e$. As $D^{\lambda_k}$ is self--dual, by taking duals we 
also have
$$[\Rad P^{\lambda_k}/\Rad^2 P^{\lambda_k}:D^{\lambda_{k+1}}]\ne0,$$
for $1\le k<e-1$. Combining these facts we conclude that 
\begin{gather*}
\Rad P^{\lambda_1}/\Rad^2P^{\lambda_1}
             \supset D^{\lambda_1}\oplus D^{\lambda_2},\\
\Rad P^{\lambda_k}/\Rad^2P^{\lambda_k}\supset 
    D^{\lambda_{k-1}}\oplus D^{\lambda_k}\oplus D^{\lambda_{k+1}},
      \text{\ for\ } 2\le k<e-1, 
\end{gather*} 
and that
$\Rad P^{\lambda_{e-1}}/\Rad^2P^{\lambda_{e-1}}
        \supset D^{\lambda_{e-1}}\oplus D^{\lambda_{e-2}}$.
Therefore, the separation diagram of~$B$~contains
$$ \psline(0,0)(0.8,0.4)(1.6,0)(1.6,0.4)(0.8,0)(0,0.4)(0,0)
   \psline(0.8,0)(0.8,0.4)
   \psline(1.6,0)(2.2,0.3)\psline(1.6,0.4)(2.2,0.1)
   \rput(2.5,0.2){\hbox{$\cdots$}}
   \psline(2.8,0.1)(3.6,0.4)(4.4,0)(4.4,0.4)(3.6,0)(2.8,0.3)
   \psline(3.6,0)(3.6,0.4)
\hspace*{42mm}$$
Consequently, $B$ has infinite representation type by
Theorem~\ref{Gabriel}.

This completes the proof of Theorem~\ref{n=e}.
\end{proof}


\section{The representation type when $n\ge2f+4$}

Recall that we are assuming that $e\ge3$ and $0\le f\le\frac e2$.
This section is devoted to the proof of the following result.

\begin{Theorem}
Suppose $n\ge2f+4$. Then $\H$ has infinite representation type.
\label{2f+4}
\end{Theorem}

\begin{proof}
By Corollary~\ref{small} it is enough to show that $\H_n$ has infinite
representation type when $n=2f+4$. Also, by the last section we may
assume that $e>n$. Let $B$ be the block $B$ with residues
$\{-1,0,0,1,1,\dots,f,f,f+1\}$; so $B$ contains the bipartition
$\((2,2),(2^f)\))$. We will show that $B$ has infinite representation
type. To do this we need to consider two cases separately.

\Case{Suppose that $f=0$ $($and $e>n=4)$} 
We consider the block~$B$ with residues $\{-1,0,0,1\}$. It is easy to see
that there are precisely six bipartitions in this block; namely,
$\lambda_1=\((0),(2^2)\)$, $\lambda_2=\((1),(2,1)\)$,
$\lambda_3=\((1^2),(2)\)$, $\lambda_4=\((2),(1^2)\)$,
$\lambda_5=\((2,1),(1)\)$ and $\lambda_6=\((2^2),(0)\)$.  Furthermore,
of these bipartitions only $\lambda_1$ and $\lambda_2$ are Kleshchev.
Since
\begin{align*}
F_0F_1F_{e-1}F_0\((0),(0)\)
   &=F_0\Big[\((0),(2,1)\)+v\((2,1),(0)\)\Big]\\
   &=\lambda_1+v\lambda_2+v\lambda_5+v^2\lambda_6\\
\intertext{and}
F_1F_{e-1}F_0^{(2)}\((0),(0)\)
   &=F_1F_{e-1}\((1),(1)\)
    =F_1\Big[\((1),(1^2)\)+v\((1^2),(1)\)\Big]\\
   &=\lambda_2+v\lambda_4+v\lambda_3+v^2\lambda_5,
\end{align*}
by Corollary~\ref{canonical} the
transpose of the decomposition matrix of $B$ is
$$\begin{array}{l|*6c} &S^{\lambda_1}& S^{\lambda_2}& S^{\lambda_3}& 
             S^{\lambda_4}& S^{\lambda_5}& S^{\lambda_6}\\\hline
D^{\lambda_1}&  1&1&0&0&1&1\\
D^{\lambda_2}&  0&1&1&1&1&0
\end{array}$$

Let $\mu_1=((0),(2,1))$ and $\mu_2=((1),(1^2))$.  We first observe
that $P^{\lambda_2}$ is the induced module $P^{\mu_2}\uparrow^B$; this
follows from a calculation in the Grothendieck group. Similarly, we
also have that $P^{\lambda_1}=P^{\mu_1}\uparrow^B$.  Further,
$P^{\mu_i}$ is uniserial of length $2$ and $[P^{\mu_i}]=2[D^{\mu_i}]$,
for $i=1,2$. Next note that $D^{\mu_i}=S^{\mu_i}$, for $i=1,2$; hence,
the branching rules for Specht modules imply that
$[D^{\mu_1}\uparrow^B]=2[D^{\lambda_1}] +[D^{\lambda_2}]$ and
$[D^{\mu_2}\uparrow^B]=[D^{\lambda_1}]+2[D^{\lambda_2}]$. 

Next note that there are surjective homomorphisms
$P^{\lambda_i}\longrightarrow D^{\mu_i}\uparrow^B$ for $i=1,2$. As
$D^{\mu_i}\uparrow^B$ is self--dual, $D^{\mu_i}\uparrow^B$
is a uniserial module whose top and bottom is isomorphic to
$D^{\lambda_i}$. To conclude, each $P^{\lambda_i}$ has a submodule
$M_i$ such that both $P^{\lambda_i}/M_i$ and $M_i$ are isomorphic to
the uniserial module $D^{\mu_i}\uparrow^B$. 
 
Now we prove that $\Ext^1(D^{\lambda_2}, D^{\lambda_2})\ne0$. As
$\Ext^1(D^{\lambda_1}, D^{\lambda_1})\ne0$, we know
that the radical series of $P^{\lambda_1}$ has the form 
$$\begin{array}{c}
D^{\lambda_1}\\
D^{\lambda_1}\oplus D^{\lambda_2}\\ 
D^{\lambda_1}\oplus D^{\lambda_2}\\ 
D^{\lambda_1}
\end{array}$$
Assume that $\Ext^1(D^{\lambda_2},D^{\lambda_2})=0$. 
Then $P^{\lambda_2}$ must have the form 
$$\begin{array}{c}
D^{\lambda_2}\\
D^{\lambda_1}\\
D^{\lambda_2}\oplus D^{\lambda_2}\\
D^{\lambda_1}\\
D^{\lambda_2}
\end{array}$$
Consider a homomorphism 
$P^{\lambda_2}\longrightarrow D^{\lambda_2}\subset 
\Rad P^{\lambda_1}/\Rad^2 P^{\lambda_1}$. 
The  assumption that 
$\Rad P^{\lambda_2}/\Rad^2 P^{\lambda_2}=D^{\lambda_1}$ together with
the radical structure of $P^{\lambda_1}$ imply that this map lifts to a
homomorphism $P^{\lambda_2}\longrightarrow P^{\lambda_1}$, whose image
is of the form 
$$
\begin{array}{c}
D^{\lambda_2}\\
D^{\lambda_1}\\
D^{\lambda_1}
\end{array}$$
This contradicts our assumption that
$\Rad^2 P^{\lambda_2}/\Rad^3 P^{\lambda_2}
          =D^{\lambda_2}\oplus D^{\lambda_2}$. 

Hence, $\Ext^1(D^{\lambda_2}, D^{\lambda_2}0$ and
$\Ext^1(D^{\lambda_1}, D^{\lambda_1})$ are both non--zero. Therefore,
$\Rad P^{\lambda_1}/\Rad^2P^{\lambda_1}$ and 
$\Rad P^{\lambda_2}/\Rad^2P^{\lambda_2}$ both contain
$D^{\lambda_1}\oplus D^{\lambda_2}$ as an $\H$--submodule. Therefore,
the separation diagram of $B$ is\,
$\psline(0,-0.1)(0.8,0.3)(0.8,-0.1)(0,0.3)(0,-0.1)$\phantom{graph}.
Consequently, $B$ is of infinite type, by Gabriel's theorem,
and~$\H$ has infinite representation type by
Lemma~\ref{subalgebra}(ii).

\Case{Suppose that $0<f\le\frac e2$} 
Now we are considering the block $B$ of $\H$ with residues
$\{-1,0,0,1,1,\dots,f,f,f+1\}$. Notice that $f+2$ is a residue in $B$
only if $f+1\ge e-1$ and $-2$ is a residue in $B$ only if $e-2\le
f+1$; however, $e>6$ because $e>n=2f+4$ and $f\ge1$; so $f+2$ and $-2$
cannot be residues in $B$. Consequently, all of the bipartitions in
$B$ have the form $\((a,b),(2^k1^l)\)$.

Define bipartitions $\lambda_1=\((0),(2^{f+2})\)$,
$\lambda_2=\((1),(2^{f+1},1)\)$, $\lambda_3=\((1^2),(2^{f+1})\)$,
$\lambda_4=\((2),(2^f,1^2)\)$ and $\lambda_5=\((2,1),(2^f,1)\)$.
We will show that $B$ has infinite representation type by computing
$\End_\H(P^{\lambda_2})$.

\begin{Lemma}
The bipartitions $\lambda_1,\lambda_2,\lambda_3,\lambda_4$
and $\lambda_5$ are all Kleshchev. Furthermore, 
$[P^{\lambda_2}]=[S^{\lambda_2}]+[S^{\lambda_3}]
                            +[S^{\lambda_4}]+[S^{\lambda_5}]$
and the first five rows of the decomposition matrix of $B$ are as
follows $($all omitted entries are zero$)$:
$$\begin{array}{l|*5c}
&D^{\lambda_1}&D^{\lambda_2}&D^{\lambda_3}&D^{\lambda_4}&D^{\lambda_5}
\\\hline
S^{\lambda_1}& 1 & . & . & . & .\\
S^{\lambda_2}& 1 & 1 & . & . & .\\
S^{\lambda_3}& 0 & 1 & 1 & . & .\\
S^{\lambda_4}& 0 & 1 & 0 & 1 & .\\
S^{\lambda_5}& 1 & 1 & 1 & 1 & 1
\end{array}$$
\end{Lemma}

\begin{proof}
We use Corollary~\ref{canonical} to compute 
$[P^{\lambda_i}]$, for $i=1,\dots,5$. We find that
$$\begin{array}{ll}
\multicolumn2c{F_1\dots F_{f+1}F_{e-1}F_0^{(2)}F_1\dots F_f\((0),(0)\)
=F_1\dots F_{f+1}F_{e-1}F_0^{(2)}\((0),(1^f)\)\qquad}\\[1mm]
\qquad&=F_1\dots F_{f+1}
           \Big[\((1),(1^{f+2})\)+v\((1^2),(1^{f+1})\)\Big]\\[1mm]
&=F_1\Big[\((1),(2^f,1^2)\)+v\((1^2),(2^f,1)\)\Big]\\[1mm]
&=\((1),(2^{f+1},1)\)+v\((2),(2^f,1^2)\)
         +v\((1^2),(2^{f+1})\)\\
      &\qquad\qquad+v^2\((2,1),(2^f,1)\).
\end{array}$$
As in Case~1, by Corollary~\ref{canonical} this shows that $\lambda_2$
is Kleshchev, proves the formula for $[P^{\lambda_2}]$ and thus gives
the second column of the decomposition matrix of $B$. The remaining
claims follow from the following calculations, which we leave to the
reader.
\begin{align*}
F_0\dots F_{f+1}F_{e-1}F_0\dots F_f\((0),(0)\)
 &=\lambda_1+v\lambda_2+v\lambda_5+\dots\\
F_1\dots F_{f+1}F_0\dots F_fF_{e-1}F_0\((0),(0)\)
 &=\lambda_3+v\lambda_5+\dots\\
F_{e-1}F_2\dots F_{f+1}F_0F_1^{(2)}F_2\dots F_fF_0\((0),(0)\)
 &=\lambda_4+v\lambda_5+\dots\\
\end{align*}
Here, ``$+\dots$'' indicates a linear combination of 
more dominant terms and if $e=\infty$ then we replace $F_{e-1}$ 
by $F_{-1}$.
\end{proof}

Consequently, 
$[P^{\lambda_2}]=4[D^{\lambda_2}]+2[D^{\lambda_1}]+2[D^{\lambda_3}]
                     +2[D^{\lambda_4}]+[D^{\lambda_5}]$.
Now, by Corollary~\ref{specht filtration}(ii), $P^{\lambda_2}$ has a
submodule isomorphic to $S^{\lambda_5}$. The composition 
factors $D^{\lambda_3}$ and $D^{\lambda_4}$ of $S^{\lambda_5}$ 
cannot appear in $\Rad P^{\lambda_2}/\Rad^2P^{\lambda_2}$ because 
then~$D^{\lambda_5}$ would appear in the head of $P^{\lambda_2}$. 
On the other hand, because $S^{\lambda_2},S^{\lambda_3}$ and 
$S^{\lambda_4}$ are indecomposable, 
$D^{\lambda_1}\oplus D^{\lambda_3}\oplus D^{\lambda_4}$ appears 
in $\Rad P^{\lambda_2}/\Rad^2P^{\lambda_2}$. Recall, again, that 
$P^{\lambda_2}/S^{\lambda_5}$ has a Specht filtration with 
successive quotients $S^{\lambda_2}, S^{\lambda_3}$ and~$S^{\lambda_4}$. 
Therefore, the $\H$--modules $D^{\lambda_3}$ and $D^{\lambda_4}$ which
appear in $\Rad P^{\lambda_2}/\Rad^2P^{\lambda_2}$ are composition
factors of $S^{\lambda_3}$ and $S^{\lambda_4}$. 
Since 
$\Rad S^{\lambda_3}=\Rad S^{\lambda_4}=D^{\lambda_2}$ and other 
$D^{\lambda_2}$ are the head and the socle of $P^{\lambda_2}$, 
$D^{\lambda_2}$ does not appear in 
$\Rad P^{\lambda_2}/\Rad^2 P^{\lambda_2}$. Thus 
$\Ext^1(D^{\lambda_2},D^{\lambda_2})=0$ and 
we see that
$P^{\lambda_2}/\Rad^3P^{\lambda_2}$ contains 
$D^{\lambda_2}\oplus D^{\lambda_2}$ as an $\H$--submodule by the same 
argument as before.  Consequently, 
$\End_\H(P^{\lambda_2}/\Rad^3 P^{\lambda_2})\cong
    k[x,y]/\<x^2,xy,y^2\>$ and~$B$ is of infinite representation type
by Lemma~\ref{end-finite}. Therefore, $\H$ has infinite representation
type by Lemma~\ref{subalgebra}(ii). 

This completes the proof of Theorem~\ref{2f+4}.
\end{proof}


\section{Finite representation type}

In this final section we show that $\H_q(B_n)$ is of finite
representation type when $n<\min\{e,2f+4\}$. To do this we use a
different combinatorial description of bipartitions which was
suggested to the first author by Fomin. Recall that we are assuming
that $q$ is a primitive $e^{\text{th}}$ root of unity in~$R$ and
that~$T_0$ satisfies the relation $(T_0-1)(T_0-q^f)=0$ where 
$0\le f\le\tfrac e2$. We note that if $K$ is a field of characteristic
zero then there is a different argument by Geck~\cite[Corollary
9.7]{Ge}. 

First, consider a partition $\lambda$. The {\sf diagram} of $\lambda$
is the set
$$\set{(i,j)\in\N^2|1\le j\le\lambda_i\And i\ge1},$$ 
which we think of as an array of boxes in the plane. Just as we can
identify $\lambda$ with its diagram we can also identify $\lambda$
with its border 
$$\set{(i,j)\in\N^2|\lambda_i<j\le\lambda_{i-1}+1}$$
(we set $\lambda_0=\infty$).  We can think of the border of $\lambda$
as a (doubly infinite) path from $(\infty,1)$ to $(1,\infty)$. Writing
$0$ for each vertical edge and $1$ for each  horizontal edge in the
border we identify $\lambda$ with a doubly infinite sequence of $0$'s
and $1$'s. We call this the {\sf path sequence} of $\lambda$.  A path
sequence is also called a {\sf Maya diagram}. 

For example, if $\lambda$ is the partition $(4,2,1)$ then by looking at
the diagram of $\lambda$ 
$$\psline(0.05,0.2)(0.6,0.2)(0.6,0.65)(1.25,0.65)(1.25,1.05)%
         (2.5,1.05)(2.5,1.45)(0.05,1.45)(0.05,0.2)
\begin{array}{*8c}
\times& \times& \times& \times& 0& 1& 1& \dots\\
\times& \times&      0&      1& 1\\
\times&      0&      1\\
     0&      1\\
     0\\
     \vdots
\end{array}$$
the path sequence of $\lambda$ is $\dots00101011011\dots$.
Here the crosses mark the nodes in the diagram of $\lambda$ and the
$0$'s and $1$'s are the nodes in the border of $\lambda$.

In order to keep track of the contents of the nodes we insert a
horizontal bar into the path sequence after the node which appears on
the diagonal $\set{(i,i)|i\in\N}$. The example above becomes
$\dots00101|011011\dots$. The bar divides the path sequence into two
regions which we refer to as the left and right regions of the path.

Note that the number of $1$'s in the left region is always the same as
the number of~$0$'s in the right region. Conversely, any sequence
$$\dots p_{-2}p_{-1}p_0|p_1p_2p_3\dots$$
of zeros and ones with $P_-=\sum_{i\le0}p_i<\infty$, $P_+=\sum_{i>0}
(1-p_i)<\infty$ and $P_-=P_+$ always corresponds to a partition.

A rim hook in (the diagram of) $\lambda$ corresponds to a subsequence
$1\dots0$ in the path sequence of $\lambda$; explicitly, if
$x=(i,j)\in[\lambda]$ then the rim hook $r_x$ corresponds to the
border path from $(\lambda_j'+1,j+1)$ to $(i,\lambda_i+1)$.  As
there is a bijection between the nodes and rim hooks of $\lambda$, the
number of subsequences of a path sequence of the form $1\dots0$ is
equal to $n$. Hence, it follows that removing a rim hook
from~$\lambda$ is the same as swapping a $0$ and a~$1$ in the path
sequence for $\lambda$.  Further, the leg length of the hook is equal
to the number of zeros in the corresponding $1\dots0$ subsequence
minus $1$. Conversely, wrapping a rim hook onto $\lambda$ is the same
as changing a~$0\dots1$ subsequence to $1\dots0$. These observations will
allow us to rephrase the Jantzen sum formula~(\ref{sum formula}) in
terms of path sequences.

Given 
$\mathbf{p}=\dots p_{-2}p_{-1}p_0|p_1p_2p_3\dots$ let 
$\bar{\mathbf p}=\dots {\bar p}_{-2}{\bar p}_{-1}{\bar p}_0|%
             {\bar p}_1{\bar p}_2{\bar p}_3\dots$ be the sequence
of left partial sums; that is, ${\bar p}_i=\sum_{j\le i}p_i$. For
example, if~$\lambda=(0)$ then the path sequence is 
$\mathbf 0=\dots 000|111\dots$ and 
$\bar{\mathbf 0}=\dots 000|123\dots$. 

The {\sf content} of a node $x=(i,j)\in\N^2$ is $c(x)=j-i$.
For a partition~$\lambda$ let
$c_k(\lambda)=\#\set{x\in[\lambda]|c(x)=k}$. Using the notation
of the last paragraph, if 
$\mathbf{p}=\dots p_{-2}p_{-1}p_0|p_1p_2p_3\dots$ is the path sequence
of $\lambda$ then the content multiplicities $c_k(\lambda)$
are given by the sequence
$\dots c_{-2}c_{-1}c_0|c_1c_2\dots=\bar{\mathbf p}-\bar{\mathbf 0}$.

We now turn to path sequences for bipartitions. If
$\lambda=(\lambda^{(1)},\lambda^{(2)})$ is a bipartition then we have
a path sequence $\dots p_{-2}p_{-1}p_0|p_1p_2\dots$ for
$\lambda^{(1)}$ and a path sequence $\dots
s_{-2}s_{-1}s_0|s_1s_2\dots$ for $\lambda^{(2)}$. Recall that the 
{\sf content} of node $x=(i,j,k)$ in $[\lambda]$ is defined to be
$c(x)=j-i+(k-1)f$ (and $\res(x)=c(x)\pmod e$).  In particular, $p_0$
corresponds to a node of content $0$ and $s_0$ corresponds to a node
of content $f$.  Accordingly, we shift the nodes in the path sequence
for $\lambda^{(2)}$ by $f$ positions to the right and define the path
sequence of the bipartition $\lambda$ to be the sequence
$\{(p_i,s_{i-f})\}$ of ordered pairs which we write with two
separation bars as follows
$$\begin{array}{*3{c@{\,}}|@{\,}*3{c@{\,}}|@{\,}*4{c@{\,}}} 
  \dots& p_{-1}  & p_0   &p_1 &\dots&  p_f& p_{f+1}& p_{f+2}&\dots\\ 
  \dots& s_{-1-f}& s_{-f}&s_{1-f}&\dots&  s_0& s_1    & s_2    &\dots 
\end{array}.$$ 
We refer the three regions in the path sequence of a bipartition
as the left, middle and right regions of the sequence.

For example, if $\lambda=\((4,2,1), (2^2,1)\)$ and $f=2$ then the 
contents in $\lambda$ and its path sequence are as follows.
$$\bitab(0&1&2&3\cr-1&0\cr -2|2&3\cr 1&2\cr 0)\quad\And\quad
  \begin{array}{r@{\,}|@{\,}c@{\,}|@{\,}l}
    \dots000101&01&1011\dots\\
    \dots000001&01&0011\dots
  \end{array}.$$

As before, it is easy to see that we can recover the contents
of a bipartition from the partial sums 
$\sum_{j\le i}(p_j+s_{j-f})$ from the path sequence.

We now develop a calculus with which to analyze path sequences of
bipartitions. Let $A=\Pair01$, $B=\Pair10$, $C=\Pair00$ and
$D=\Pair11$  be the four possible ordered pairs which can appear in
the path sequence of a bipartition. Define $a_l$, $a_m$ and $a_r$ to
be the number of $A$'s in the left, middle and right regions,
respectively, of the path sequence; similarly, we define
$b_l,b_m,b_r,c_l,c_m,c_r,d_l,d_m$ and $d_r$. Notice that~$c_l$
and~$d_r$ are both infinite; all of the other quantities are
non--negative and finite.

\begin{Point}
*\label{inequalities}
Suppose that $n\le 2f+3$. Then
\begin{enumerate}
\item $b_l+d_l=a_m+c_m+a_r+c_r$.
\item $b_r+c_r=a_l+d_l+a_m+d_m$.
\item $f=a_m+b_m+c_m+d_m$.
\item $2f+3\ge (b_l+d_l)(a_m+c_m+a_r+c_r) 
              +(a_l+d_l+a_m+d_m)(b_r+c_r)
              \\\phantom{2f+3\ge}+(b_m+d_m)(a_r+c_r)
              +(a_l+d_l)(b_m+c_m)
              +a_lb_l+a_mb_m+a_rb_r$.
\end{enumerate}
\end{Point}

\begin{proof}
Parts (i) and (ii) both follow from the fact mentioned earlier
that the number of ones in a region to the left of a bar is equal to
the number of zeros in the region to the right of the same bar. Part
(iii) is true by definition. Finally, (iv) follows by counting the rim
hooks in $\lambda$ (recall that the rim hooks in the path sequence for
a partition correspond to subsequences of the form $1\dots0$); this is
the only place where we use the restriction on $n$.
\end{proof}

We want to understand this system of inequalities when $n\le 2f+3$. As
a first step, adding parts (i) and (ii) and subtracting (iii) we see
that
$$b_l+b_m+b_r=f+a_l+a_m+a_r.\Number{b=f+a}$$
%
Careful inspection shows that by combining parts (i),
(ii) and (iv) of~(\ref{inequalities}), and omitting some terms for the
second inequality, we have
\begin{align*}
2f+3&\ge(a_m+c_m+a_r+c_r)^2+(a_l+d_l+a_m+d_m)^2+(b_m+d_m)a_r\\
    &\qquad+(b_m+d_m)c_r+(a_l+d_l)(b_m+c_m)
        +a_lb_l+a_mb_m+a_rb_r\\
    &\ge(a_l+a_m+a_r)(a_m+b_m+c_m+d_m)+a_la_m+a_ma_r+a_l^2+a_m^2\\
    &\qquad+a_r^2+a_l(b_l+d_l)+a_r(b_r+c_r)\\
    &=(a_l+a_m+a_r)f+a_la_m+a_ma_r+a_l^2+a_m^2+a_r^2\\&\qquad
     +a_l(a_m+c_m+a_r+c_r)+a_r(a_l+d_l+a_m+d_m);
\intertext{the last line again uses (i)--(iii) of 
(\ref{inequalities}). Finally, throwing away a few more terms and
rearranging gives}
2f+3&\ge(a_l+a_m+a_r)(f+a_l+a_m+a_r).
\end{align*}
Now, if $a_l+a_m+a_r\ge2$ then the right hand side is greater than or
equal to $2f+4$; as this is impossible we must therefore have 
$$\lefteqn{a_l+a_m+a_r\le 1.}\Number{small a}$$
In particular, note that at most one of $a_l$, $a_m$ and $a_r$ can be
non--zero. Notice that this is the first time that the assumption
$n\le2f+3$ has really been needed; it is exactly what is required to
ensure that $a_l+a_m+a_r\le 1$.

Using similar arguments it is possible to classify the possible path
sequences when $n\le 2f+3$; however, we won't need this.

\begin{Theorem}
Suppose that $n<\min\{e,2f+4\}$. Then $\H$ is of finite
representation type.
\end{Theorem}

\begin{proof}
Fix a block $B$ of $\H$ and recall from Corollary~\ref{blocks} that two
simple modules belong to the same block only if they have a common 
multiset of residues.  Note that because $n<e$ the number of distinct
residues contained in the diagram of $\lambda$ is strictly less
than~$n$; consequently, we can find a~$k$, with $0\le k<e$, such that
$-k-1\pmod e$ is not a residue in $B$.

Suppose that $D^\lambda$ appears in $B$. Then the contents of all of
the nodes in $\lambda^{(1)}$ are contained in the interval
$[-k,e-k-1]$.

\Case{Suppose that $e-k\ge f$} As $-k-1\pmod e$ is not a residue for
the block~$B$ the contents of the nodes in $\lambda^{(2)}$ are all
contained in the interval $[-k,e-k]$~--- note that $f\in[-k,e-k]$.
Therefore, the multiset of residues for the block $B$ is the same as
the multiset of contents for $B$. Consequently, we can unwrap a rim
hook from~$\lambda$ and wrap it back on again without changing the
residue of the foot node only if the resulting bipartition $\mu$ has
the same multiset of contents as $\lambda$. 

Recall that unwrapping a rim hook from a partition is the same as
swapping the ends of a $1\dots0$ subsequence to give $0\dots1$ and
that wrapping a hook back on changes some $0\dots1$ into $1\dots0$.
Now, the contents of a bipartition $\lambda$ are determined by the
partial sums in the path sequence of $\lambda$; because of this, the
only way to unwrap a  rim hook from $\lambda$ and then wrap it back on
to give a bipartition $\mu$ with the same multiset of contents is by
interchanging some $A$ and $B$ in the path sequence:
$$\lambda=\ldots B\ldots A\ldots 
	\qquad\To\qquad \mu=\ldots A\ldots B\ldots.$$ 
Moreover, $\lambda\gdom\mu$ if and only if $A$ moves to the left.
(Note that $|\mu|=|\lambda|$ in this case as the number of $1\dots0$
subsequences in the two path sequences is the same.)

If $a_l+a_m+a_r=0$ then the path sequence for $\lambda$ does not
contain any $A$'s so by the sum formula, Theorem~\ref{sum formula}, 
and by Proposition~\ref{cor of sum formula}, 
the Specht module $S^\lambda=D^\lambda=P^\lambda$ is the only simple
module in the block $B$. In particular, $B$ is semisimple and so of
finite type in this case.

If $a_l+a_m+a_r\ne0$ then $a_l+a_m+a_r=1$ by (\ref{small a}).  In this
case by Proposition~\ref{cor of sum formula} and (\ref{b=f+a}) the block~$B$
contains at most $f+2$ bipartitions; namely, the bipartitions
$\lambda_0,\dots,\lambda_{f+1}$ whose path sequences contain exactly
one $A$ and $(f+1)$~$B$'s and which agree with the path sequence for
$\lambda$ on all of the $C$'s and $D$'s. By ordering these
bipartitions according to the location of the (unique) $A$ in their
path sequence we may assume that
$\lambda_{f+1}\gdom\dots\gdom\lambda_0$ (for example, $A$ occupies the
leftmost position in the path sequence of $\lambda_0$ and the
rightmost position in $\lambda_{f+1}$).  

For our purposes, it is enough to prove that
$S^{\lambda_0}=D^{\lambda_0}$ and
$[S^{\lambda_i}]=[D^{\lambda_i}]+[D^{\lambda_{i-1}}]$ for $0<i\le f$;
in particular, we do not need to know that $D^{\lambda_i}\ne0$, for
$0\le i\le f$. In fact these modules are always non--zero; we include
the proof below because it yields the remarkable fact that when
$n<\min\{e,2f+4\}$ the number of Specht modules belonging to a block
is either $1$, $f+2$ or $e-f+2$.

The removable nodes in a partition correspond to the $10$
subsequences in the path sequence. Suppose that $0\le i\le f$. Then
$\lambda_i$ contains a removable node~$x\in[\lambda_i^{(2)}]$;
furthermore, this node is automatically good because if $x$ is a
$r$--node then there is no addable $r$--node below $x$ because
$|\lambda^{(2)}|<e$.  Let~$\mu_i$ be the
bipartition with $[\mu_i]=[\lambda_i]\setminus\{x\}$. Then~$\mu_i$ is
Kleshchev because either the path sequence for $\mu_i$ contains an
$A$ (so $\mu_i$ is Kleshchev by induction on~$n$), 
or~$D^{\mu_i}=S^{\mu_i}\ne0$ (by the second last paragraph); hence,
$\lambda_i$ is Kleshchev. On the other
hand, $\lambda_{f+1}$ is not Kleshchev because either we can apply
induction after removing a node from~$\lambda_{f+1}^{(1)}$, or the
path sequence for~$\lambda_{f+1}$ is $\dots BBBA\dots$ in which case
it is easy to see that~$\lambda_{f+1}$ is not Kleshchev.

Now $S^{\lambda_0}=D^{\lambda_0}$ by Theorem~\ref{DJM}(iii)
since $\lambda_i\gdom\lambda_0$ for $i>0$.
To complete the proof we claim that 
$[S^{\lambda_i}]=[D^{\lambda_i}]+[D^{\lambda_{i-1}}]$,
for~$i=1,\dots,f$, and $S^{\lambda_{f+1}}=D^{\lambda_f}$. To see
this we apply the sum formula. As discussed earlier, the leg length of
a hook $1\ldots 0$ in a path sequence is given by the number of $0$'s
strictly contained in the subsequence.  Consequently, when we unwrap
the hook $B\ldots A$ from $\lambda_i$ and wrap it back on to give some
$\lambda_l$ then, modulo~$2$, the difference in the leg lengths of
the two rim hooks is equal to the number of $B$'s which are
strictly contained in the subsequence for the rim hook. Therefore, by
Theorem~\ref{sum formula}, for $i=1,\dots,f+1$
$$\sum_{j>0}[S^{\lambda_i}(j)]=[S^{\lambda_{i-1}}]-[S^{\lambda_{i-2}}]
                    +\dots+(-1)^{i-1}[S^{\lambda_0}].$$
As we already know that $\lambda_0,\dots\lambda_f$ are Kleshchev, and
that $\lambda_{f+1}$ is not, our claim now follows by
induction on $i$. Consequently, the decomposition matrix of the 
block~$B$ is 
$$\begin{array}{l|cccc}
&D^{\lambda_0} & D^{\lambda_1}&\dots&D^{\lambda_f}\\\hline
S^{\lambda_0}     & 1 & \\ 
S^{\lambda_1}     & 1 & 1\\ 
\vdots            &   &\ddots&\ddots\\ 
S^{\lambda_f}     &   &   & 1 & 1\\ 
S^{\lambda_{f+1}} &   &   &   & 1 
\end{array}$$ 
and $B$ has finite representation type by Theorem~\ref{weight one}.

\Case{Suppose that $0<e-k<f$} 
In this case the contents of the nodes in~$\lambda^{(2)}$ are
contained in the interval $[e-k,2e-k-1]$. Renormalizing $T_0$ as
$q^fT_0$ the relation for~$T_0$ becomes $(T_0-1)(T_0-q^{e-f})=0$. The
Specht module $S^\lambda$ is relabelled
as~$S^{(\lambda^{(2)},\lambda^{(1)})}$ and the residues in $[\lambda]$
are all changed by adding $e-f\equiv-f\pmod e$. Consequently, the
residues for $\lambda$ are all contained in the interval
$[e-f-k,2e-f-k-1]$.  Therefore, with this renormalization, the
multiset of residues for $B$ is the same as the multiset of contents
for $B$.  Consequently, we can repeat the argument of Case~1 to deduce
that decomposition matrix for $B$ has the form above; so, again,~$B$
has finite representation type by Theorem~\ref{weight one}.  
\end{proof}

\let\em\it


\end{document}